\documentclass[11pt]{amsart}
\usepackage{graphicx}
\usepackage{amsfonts, amsmath, amsfonts, amssymb}
\newtheorem{thm}{Theorem}[section]

\theoremstyle{definition}
\newtheorem{defn}[thm]{Definition}
\theoremstyle{remark}
\newtheorem{rem}[thm]{Remark}
\numberwithin{equation}{section}

\numberwithin{equation}{section}

\numberwithin{equation}{section} \theoremstyle{cond}

 \numberwithin{equation}{section}
 \newtheorem{quest}[thm]{Question}

\begin{document}


\title[]{The existence of thick triangulations -- an ``elementary'' proof}%
\author{Emil Saucan \and Meir Katchalski}%
\address{Mathematics Department, Technion, Technion City, Haifa 32000 }%
\email{semil@tx.technion.ac.il}%
\email{meirk@tx.technion.ac.il}%

\subjclass{AMS Classification. 53A07, 57R05, 57M60, 58A05}
\keywords{thick triangulation,  $\mathcal{C}^1$ manifolds, cut locus}%

\date{\today}%
\begin{abstract}
We provide an alternative, simpler  proof of the existence of thick
triangulations for noncompact $\mathcal{C}^1$ manifolds. Moreover,
this proof is simpler than the original one given in \cite{pe},
since it mainly uses tools of elementary differential topology. The
role played by curvatures in this construction is also emphasized.
\end{abstract}

\maketitle


\section{Introduction}

The existence of the so called ``thick'' or ``fat'' triangulations
is an important both in Pure Mathematics, in Differential Geometry
(where it plays a crucial role in the computation of curvatures for
piecewise-flat approximations of smooth Riemannian manifolds, with
applications to Regge Calculus, see \cite{reg:61}); and in Geometric
Function Theory (mainly in the construction of {\it quasimeromorphic
mappings}, see \cite{ms2}, \cite{tu}, \cite{pe}, \cite{s1},
\cite{s4}, \cite{s3}).
Thick triangulations play also an important role in many
applications, mainly in Computational Geometry, Computer Graphics,
Image Processing and their related fields (see, \cite{ab}, \cite{bcer}, \cite{cdr}, \cite{e}, \cite{pa}, \cite{saz}, to name just a few). 

Recall that ``thick'' or ``fat'' triangulations are defined as
follows:

\begin{defn} Let $\tau \subset \mathbb{R}^n$ ; $0 \leq k \leq n$ be a $k$-dimensional simplex.
The {\it thickness}  $\varphi$ of $\tau$ is defined as being:
\begin{equation}
\varphi = \varphi(\tau) = \hspace{-0.3cm}\inf_{\hspace{0.4cm}\sigma
< \tau
\raisebox{-0.25cm}{\hspace{-0.9cm}\mbox{\scriptsize$dim\,\sigma=j$}}}\!\!\frac{Vol_j(\sigma)}{diam^{j}\,\sigma}
\end{equation}
The infimum is taken over all the faces of $\tau$, $\sigma < \tau$,
and $Vol_{j}(\sigma)$ and $diam\,\sigma$ stand for the Euclidian
$j$-volume and the diameter of $\sigma$ respectively. (If
$dim\,\sigma = 0$, then $Vol_{j}(\sigma) = 1$, by convention.)
\\ A simplex $\tau$ is $\varphi_0${\it-thick}, for some $\varphi_0 > 0$, if $\varphi(\tau) \geq \varphi_0$. A triangulation (of a submanifold of $\mathbb{R}^n$) $\mathcal{T} = \{ \sigma_i \}_{i\in \bf I}$ is
$\varphi_0${\it-thick} if all its simplices are $\varphi_0$-thick. A
triangulation $\mathcal{T} = \{ \sigma_i \}_{i\in \bf I }$ is {\it
thick} if there exists $\varphi_0 \geq 0$ such that all its
simplices are $\varphi_0${\it-thick}.
\end{defn}

The definition above is the one introduced in \cite{cms}. For some
different, yet equivalent definitions of thickness, see \cite{ca1},
\cite{ca2}, \cite{mun}, \cite{pe}, \cite{tu}.


For $\mathcal{C}^\infty$ Riemannian manifolds without boundary the
the existence has been proved by Peltonen \cite{pe}. This result was
extended by the first author to include manifolds of lower
differentiability class with boundary \cite{s2} and to a large class
of orbifolds \cite{s3}.

 Peltonen's proof is based on the construction of an exhaustion of
the manifold using a delicate curvature-based argument, both to
decide the ``size'' of the compact ``pieces'' (i.e. of the elements
of the exhaustion) and to chose the mash of the triangulation (see
\cite{pe} and also Section 2). The technique used in \cite{s2} is
much more elementary, using only the Differential Topology apparatus
(and results) of \cite{mun}. This discrepancy between methods
creates a kind of ``esthetic asymmetry'' that naturally gives rise
to the following question:

\begin{quest}
Is it possible to prove the existence of thick triangulations for
non-compact manifolds using only techniques of Elementary
Differential Topology?
\end{quest}

We shall prove that the answer to the question above is positive by
showing that the ``meshing'' technique of thick triangulations
developed in \cite{s2} allows us to discard the curvature
considerations of the original proof.

However, in discarding the curvature-related information, one also
loses geometric intuition and, with it, any possibility of applying
the technique in any non-trivial, concrete case. Hence, the next
question ensues immediately:

\begin{quest}
Can one recover the Differential Geometric information (i.e.
curvatures) from the constructed $PL$ triangulation?
\end{quest}

We show that, again, the answer is affirmative, and it follows from
the results of \cite{cms}. Moreover, we indicate how this approach
can also be simplified using tools that may be considered more
``elementary''.

 The reminder of the paper is organized as follows: In the
next Section we briefly sketch, for the benefit of the reader and
for the sake of the paper's self-containment, the main steps in the
proofs of the main results in \cite{pe} and \cite{s2}. In Section 3,
we show how our result in \cite{s2} allows us to give a simpler,
``elementary'' proof of Peltonen's theorem. Finally, in the last
Section, we discuss the role played by curvature in our simplified
proof and indicate how, using our method, one still can recapture
the Differential Geometric information encoded in Peltonen's proof.


\section{Background}

We bring below a very brief sketch of the methods used in proving
the existence of thick triangulations. We concentrate mainly on
those aspects that are pertinent to our present study.

\subsection{Open Riemannian  Manifolds}

First, a number of necessary definitions:


Let $M^n$ denote an $n$-dimensional complete Riemannian manifold,
and let $M^n$ be isometrically embedded into $\mathbb{R}^\nu$
(``$\nu$''-s existence is guaranteed by Nash's Theorem -- see, e.g.
\cite{pe}).

Let $\mathbb{B}^\nu(x,r) = \{y \in \mathbb{R}^\nu\,|\, d_{eucl} <
r\}$; $\partial\mathbb{B}^\nu(x,r) = \mathbb{S}^{\nu-1}(x,r)$. If $x
\in M^n$, let $\sigma^n(x,r) = M^n \cap \mathbb{B}^\nu(x,r)$,
$\beta^n(x,r) = exp_x\big(\mathbb{B}^n(0,r)\big)$, where: $exp_x$
denotes the exponential map: $exp_x:T_x(M^n) \rightarrow M^n$, and
where $\mathbb{B}^n(0,r) \subset T_x\big(M^n\big)$,
$\mathbb{B}^n(0,r) = \{y \in \mathbb{R}^n\,|\,d_{eucl}(y,0) < r\}$.

The following definitions generalize in a straightforward manner classical ones used for surfaces in
$\mathbb{R}^3$:

\begin{defn}
\begin{enumerate}
\item $\mathbb{S}^{\nu-1}(x,r)$ is {\em tangent} to $M^n$ at $x\in M^n$ iff there exists $\mathbb{S}^n(x,r) \subset
\mathbb{S}^{\nu-1}(x,r)$, such that $T_x(\mathbb{S}^n(x,r)) \equiv
T_x(M^n)$.
\item Let $l \subset \mathbb{R}^\nu$ be a line, then $l$ is {\em secant} to $X \subset M^n$ iff $|\,l \cap X| \geq 2$.
\end{enumerate}
\end{defn}

\begin{defn}
\begin{enumerate}
\item $\mathbb{S}^{\nu-1}(x,\rho)$ is an {\rm osculatory sphere} at $x \in M^n$ iff:
\begin{enumerate}
\item $\mathbb{S}^{\nu-1}(x,\rho)$ is tangent at x;
\\ and
\item $\mathbb{B}^n(x,\rho) \cap M^n = \emptyset$.
\end{enumerate}
\item Let $X \subset M^n$. The number $\omega = \omega_X = \sup\{\rho > 0\,|\, \mathbb{S}^{\nu-1}(x,\rho) \; {\rm osculatory} \\{\rm at\; any}\; x \in
X\}$ is called the {\em maximal osculatory} ({\em tubular}) {\em radius} at $X$.
\end{enumerate}
\end{defn}

\begin{rem}
There exists an osculatory sphere at any point of $M^n$ (see
\cite{ca3}\,).
\end{rem}

\begin{defn} Let $U \subset M^n, U \neq \emptyset$, be a relatively compact set, and let $T = \bigcup_{x \in
\bar{U}}\sigma(x,\omega_U)$. The number $\kappa_U =
\max\{r\,|\,\sigma^n(x,r)\;  {\rm  is\; connected \; for \; all}\; s
\leq \omega_U,\, x \in \bar{T}\}$, is called the {\em maximal
connectivity radius} at U, defined as follows:
\end{defn}

Note that the maximal connectivity radius and the maximal osculatory
radius are interconnected by the following inequality (\cite{pe},
Lemma 3.1):

\[\omega_U \leq \frac{\sqrt{3}}{3}\kappa_U\,.\] \label{ec:1}

We are now able to present the main steps of Peltone's proof, which
generalizes both the result and method of proof of Cairns
\cite{ca3}:

%
%
%
%
%
%
%

%
%



\begin{enumerate}
\item
\begin{enumerate} 
\item Construct an exhaustive set $\{E_i\}$ of $M^n$,
generated by the pair $(U_i,\eta_i)$, where:
\begin{enumerate}
\item $U_i$ is the relatively compact set $E_i \setminus \bar{E}_{i-1}$ and
\item $\eta_i$ is a number that controls the fatness of the simplices of the triangulation of $E_i$\,,  constructed in Part 2, such that it will not differ to much
on adjacent simplices, i.e.:
\\ (${\rm ii_1}$) The sequence $(\eta_i)_{i\geq1}$ descends to $0$\,;
\\ (${\rm ii_2}$) $2\eta_i \geq \eta_{i-1} \,.$
\end{enumerate}

The numbers $\eta_i$ are chosen such that they satisfy the following bounds:
\[\eta_i \leq \frac{1}{4}\min_{i \geq 1}\{\omega_{\bar{U}_{i-1}},\omega_{\bar{U}_i},\omega_{\bar{U}_{i+1}}\}\,.\]
\item 
\begin{enumerate}
\item Produce a maximal set $A$, $|A| \leq \aleph_0$, s.t. $A \cap U_i$ satisfies:
\\ (${\rm i_1}$) a density condition, namely:
\[d(a,b) \geq \eta_i/2\,, {\rm for\; all}\; i \geq 1\,;\]
(${\rm i_2}$) a ``gluing'' condition for $U_i, U_{i+1}$\,, i.e.
their intersection is large enough.

Note that according to the density condition (${\rm i_1}$), the
following holds: For any $i$ and for any $x \in \bar{U}_i$, there
exists $a \in A$ such that $d(x,a) \leq \eta_i/2$\,.

\item Prove that the Dirichlet complex $\{\bar{\gamma}_i\}$ defined by the sets $A_i$ is a cell complex and
every cell has a finite number of faces (so that it can be triangulated in a standard manner).
\end{enumerate}
\end{enumerate}

\item 
Consider first the dual complex $\Gamma$, and prove that it is a Euclidian
simplicial complex with a ``good'' density. Project then $\Gamma$ on
$M^n$ (using the normal map). Finally, prove that the resulting
complex $\widetilde{\Gamma}$ can be triangulated by fat simplices.
%
\end{enumerate}

\subsection{Manifolds With Boundary}
%
%
The idea of the proof in this case is to build first two fat
triangulations: $\mathcal{T}_{1}$ of a product neighbourhood $N$ of
$\partial M^n$ in $M^n$ and $\mathcal{T}_{2}$ of $int\, M^n$ (its
existence follows from Peltonen's result), and then to ``mash'' the
two triangulations into a new triangulation $\mathcal{T}$, while
retaining their thickness (see \cite{s2}). While the mashing
procedure of the two triangulations is basically the classical one
of \cite{mun}, the triangulation of $\mathcal{T}_{1}$ has been
modified, in order to ensure the thickness of the simplices of
$\mathcal{T}_{1}$ (see \cite{s2}, Theorem 2.9). To thicken
triangulations one can use either the method used in \cite{cms}, or,
alternatively, the one developed in \cite{s1}. For the technical
details, see \cite{s2}.

%


\section{Main Result}

The idea of the proof 
is to use the basic fact that $M^n$ is $\sigma$-compact, i.e. it
admits an exhaustion by compact submanifolds $\{K_j\}_j$ (see, e.g.
\cite{spi}). This is a standard fact for metrizable manifolds.
However, it is conceivable that the ``cutting surfaces'' $N_{ij}$\,,
$\bigcup_{\scriptscriptstyle{i=1,...k_j}}N_{ij} =
\partial K_j$\,, are merely $\mathcal{C}^0$, so even the existence of a
triangulation for these hypersurfaces is not always assured, 
hence a fortiori that of smooth triangulations.
(See. e.g. \cite{th} for a brief review of the results regarding the
existence of triangulations).

To show that one can obtain (by ``cutting along'') smooth
hypersurfaces, we briefly review the main idea of the proof of the
$\sigma$-compactness of $M^n$ (for the full details, see, for
example \cite{spi}): Starting from an arbitrary base point $x_0 \in
M^n$, one considers the interval $I = I(x_0) = \{r > 0\,|\,
\beta^n(x_0,r)\; {\rm is \; compact}\}$, where $\beta^n(x_0,r)$ is
as in Section 2. If $I = \mathbb{R}$, then $M^n =
\bigcup_{\scriptscriptstyle
1}^{\scriptscriptstyle\infty}{\beta^n(x,n)}$, thence
$\sigma$-compact. If $I \neq \mathbb{R}$, one constructs the
compacts sets $K_j$, $K_0 = \{x_0\}$, $K_{n+1} =
\bigcup_{\scriptscriptstyle y \in K_j}\beta^n(y,r(y))$, where $r(y)
= \frac{1}{2}\sup\{r \in I(y)\}$. Then it can be shown that $M^n =
\bigcup_{\scriptscriptstyle n \geq 0}K_j$, i.e. $M^n$ is
$\sigma$-compact.

The smoothness of the surfaces $N_{ij}$ now follows from Wolter's
result \cite{wo} regarding the $2$-differentiability of the cut
locus of the exponential map.

By \cite{ca1} (see also \cite{ca2}) the sets $K_j$ and $N_{ij}$
admit thick triangulations and, moreover, these triangulations have
thickness $\varphi_1 = \varphi_1(n)$ and $\varphi_2 =
\varphi_2(n-1)$, respectively. One can thus apply repeatedly the
``mashing'' technique developed in \cite{s2}, for collars of
$N_{ij}$ in $K_j$ and $K_{j + 1}$, $j \geq 0$, rendering a
triangulation of $M^n$, of uniform thickness $\varphi = \varphi(n)$
(see \cite{s2}, \cite{cms}).


\begin{rem}
Instead of the less known and more complicated method of \cite{ca1},
we could have employed the simpler and more intuitive one in
\cite{ca2}. However, the later one still makes use of the principal
curvatures (at each point of) the manifold, thus using this method
the ``Differential Geometric content'', so to speak, of our proof
would still be rather higher. Since we strive to obtain a proof that
is as elementary as possible, i.e. using only (or mainly) the tools
of elementary differential topology, we prefer to adopt the methods
of Cairns' earlier work.
\end{rem}

\begin{rem}
In some special cases a ``purely'' Euclidean construction can be
achieved, without resorting to embeddings and piecewise-flat (or
just piecewise-linear) approximations. Such a construction is
provided in \cite{ep}, where a method of dividing a non-compact
hyperbolic manifold into (canonical) Euclidean pieces is described.
Note that, moreover, these pieces can be easily subdivided into
thick simplices using a number of the methods mentioned above (and
in the bibliography).
\end{rem}

\begin{rem}
The opposite problem, that of extending a thick triangulation from
the interior of a manifold to its boundary, is also worth
considering, in itself and for its importance in Geometric Function
Theory (see \cite{s3}).
\end{rem}


\section{Curvatures}

We conclude with a few brief remarks regarding the role of
curvature:, its possible ``recovery'' in the $PL$ context, and
further directions of study.

First, let us note that, manifestly enough, our simplified proof is
not ``Differential Geometry free'', since we used both the
exponential map and Wolter's result. Regarding an ``elementary''
proof, this is evidently a weakness, as is the appeal we have made
to Nash's Embedding Theorem (with its complicated Differential
Geometric apparatus).
One would be tempted to substitute for the balls $\beta^n(x,j)$,
Euclidean balls $\mathbb{B}^\nu(x,j)$, and replace, if necessary,
the surfaces $L_{jk} = \partial\mathbb{B}^\nu(x,j) \cap M^n$ by
$\mathcal{C}^1$ surfaces $L^\ast_j$, that are $\varepsilon$-close to
$L_j$\,. Such an approach would allow us to consider only
$\mathcal{C}^1$ embeddings, thus permitting us to dispense with the
use of Nash's Theorem. Unfortunately, in general, one cannot discard
the use of
the exponential map (and, consequently, neither can one ignore Wolter's work). %
Indeed, for wildly embedded manifolds, the simple method above is
not applicable.
%
However, if only tame embeddings are considered, than the method
above can be employed.

Secondly, let us note that Wolter's method is -- obviously -- not
independent of curvature. Therefore, the curvature considerations do
play a decisive role, even if in a more ``soft'' manner. Moreover,
even if applying the very simple and direct method described in the
previous paragraph, the curvature plays an intrinsic role in
determining the meshes of the triangulations of two adjacent
``pieces'' $K_j$ and $K_{j+1}$, ($K_j = \mathbb{B}^\nu(x,j) \cap
M^n$), and their common boundary $L^\ast_{jk}$. For this reason, a
``naive'' exhaustion of the manifold using Euclidean balls, hence
without control of curvature, may prove, in general, to be
counterproductive, rendering not a monotone (and even highly
oscillating) series of meshes for the members of the exhaustion.

Finally, as noted in the Introduction, the (Lipschitz-Killing)
curvatures have ``good'' convergence properties, under
piecewise-flat (secant) approximations of $M^n$. This result of
Cheeger et al. (see \cite{cms}) is the goal for which they developed
the ``mashing of triangulations'' technique mentioned above. 
However, a simpler and more direct approach to curvatures
computation in $PL$ approximation, along the lines indicated in
\cite{s4} and based on metric curvatures, is currently in
preparation.


\subsection*{Acknowledgment}

The first author wishes to express his gratitude to Professor Shahar
Mendelson -- his warm support is gratefully acknowledged. He would also like to thank 
Professor Klaus-Dieter Semmler for bringing to his attention
Wolter's work.



\end{document}